\documentclass{amsart}
\usepackage{amssymb, amsmath, amsfonts, amscd}

\input xypic

\theoremstyle{plain}
\newtheorem{theorem}{Theorem}[section]

\newtheorem{lemma}[theorem]{Lemma}
\newtheorem{proposition}[theorem]{Proposition}

\newtheorem{example}[theorem]{Example}

\theoremstyle{definition}
\newtheorem{definition}[theorem]{Definition}

\theoremstyle{remark}

\numberwithin{equation}{theorem}

\renewcommand{\P}{\mathbf{P} }

\renewcommand{\H}{\operatorname{H} }

\newcommand{\K}{\operatorname{K} }

\renewcommand{\K}{\operatorname{K} }

\renewcommand{\c}{\mathbf{C} }
\renewcommand{\r}{\mathbf{R}} 
\newcommand{\Z}{\mathbf{Z}} 
\newcommand{\Q}{\mathbf{Q}} 
\begin{document}

\title{Notes on the Chern-character} 
\author{Helge Maakestad }
\address{NTNU, Trondheim} 
\email{Helge.Maakestad@math.ntnu.no}

\date{November 2006} 

\begin{abstract} Notes for some talks given at the seminar on characteristic
  classes at NTNU
  in autumn 2006. In the note a proof of the existence of a Chern-character from
  complex K-theory to any cohomology theory with values in
  graded $\Q$-algebras equipped with a theory of characteristic
  classes is given. It respects the Adams and Steenrod operations. 
\end{abstract}

\maketitle
\tableofcontents

\section{Introduction}

The aim of this note is to give an axiomatic
and elementary treatment of Chern-characters of vectorbundles with values in a class of cohomology-theories
arising in topology and algebra. Given a theory of Chern-classes for
complex vectorbundles with values in singular cohomology  one gets in a natural way
a Chern-character from complex K-theory to singular cohomology using
the projective bundle theorem  and the Newton polynomials. The
Chern-classes of a complex vectorbundle may be defined using the
notion of an Euler class (see section 14 in \cite{mil}) and one may prove that a theory of Chern-classes with
values in singular cohomology is unique. In this note it is shown
one may relax the conditions on the theory for Chern-classes and
still get a Chern-character. Hence the Chern-character depends on some choices.

Many cohomology theories which associate to a space a graded
commutative $\Q$-algebra $\H^*$ satisfy the projective bundle
property for complex vectorbundles. This is true for De Rham-cohomology of a real compact
manifold, singular cohomology of a compact topological space and
complex K-theory. The main aim of this note is to give a self contained
and elementary proof of the fact that any such
cohomology theory will recieve a Chern-character from complex K-theory
respecting the Adams and Steenrod operations. 

Complex K-theory for a topological space $B$ is considered, and
characteristic classes in K-theory and 
operations on K-theory such as the Adams operations are constructed explicitly, following \cite{karoubi}. 

The main result of the note is the following (Theorem \ref{main}):

\begin{theorem} Let $\H^*$ be any rational cohomology theory
  satisfying the projective bundle property. 
There is
for all $k\geq 1$ a commutative diagram
\[ 
\diagram  K^*_\c(B) \rto^{Ch}\dto^{\psi^k} & \H^{even}(B)
\dto^{\psi^k_H} \\
K^*_\c(B) \rto^{Ch}  & \H^{even}(B) 
\enddiagram 
\]
where $Ch$ is the Chern-character for $\H^*$, $\psi^k$ is the Adams
operation and $\psi^k_H$ is the Steenrod operation. 
\end{theorem}
The proof of the result is analogous to the proof of
existence of the Chern-character for singular cohomology.

\section{Euler classes and  characteristic classes}
In this section we consider axioms ensuring that any cohomology theory
$\H^*$ satisfying these axioms, recieve a Chern-character for
complex vectorbundles. By a cohomology theory we mean a contravariant
functor
\[ \H^*: \underline{Top}\rightarrow \underline{\Q-algebras} \]
from the category of topological spaces to the category of graded
commutative $\Q$-algebras with respect to continuous maps of
topological spaces. We say the theory satisfy the \emph{projective
  bundle property} if the following axioms are satisfied: For any
rank $n$ complex continuous vectorbundle $E$ over a compact space $B$
There is an \emph{Euler class} 
\begin{align}
u_E\in \H^2(\P(E))
\end{align}
where $\pi:\P(E)\rightarrow B$ is the projective bundle associated to $E$.
This assignment satisfy the 
following properties: The Euler class is natural, i.e for any
map of topological spaces $f:B'\rightarrow B$ it follows
\begin{align}
&f^*u_E=u_{f^*E}
\end{align}
For $E=\oplus_{i=1}^n L_i$ where $L_i$ are linebundles there is an equation
\begin{align}
& \prod_{i=1}^n (u_E-\pi^*u_{L_i})=0 \text{ in } \H^{2n}(\P(E))
\end{align}
The map $\pi^*$ induce an injection $\pi^*:\H^*(B)\rightarrow \H^*(\P(E))$ and there is an equality
\[ \H^*(\P(E))=\H^*(B)\{1,u_E,u_E^2,..,u_E^{n-1}\} .\]

Assume $\H^*$ satisfy the projective bundle property.
There is by definition an equation
\[ u^n_E-c_1(E)u_E^{n-1}+\cdots +(-1)^nc_n(E)=0 \]
in $\H^*(\P(E))$ .
\begin{definition}
The class $c_i(E)\in \H^{2i}(B)$ is the \emph{i'th characteristic class} of
$E$.
\end{definition}

\begin{example}If $\P(E)\rightarrow B$ is the projective bundle of a
  complex vector bundle and $u_E=e(\lambda(E)) \in \H^{2}(\P(E),\Z)$ is the Euler
  classe of the tautological linebundle $\lambda(E)$ on $\P(E)$ in singular
  cohomology as defined in Section 14 \cite{mil}, one verifies the properties above are satisfied. 
One gets the Chern-classes $c_i(E)\in \H^{2i}(B,\Z)$ in singular cohomology.
\end{example}

\begin{definition}\label{char}
A \emph{theory of characteristic classes} with values in a cohomology
theory $\H^*$ is an assignment
\[ E\rightarrow c_i(E)\in \H^{2i}(B) \]
for every complex finite rank vectorbundle $E$ on $B$ satisfying the
following
axioms:
\begin{align}
&\label{one}f^*c_i(E)=c_i(f^*E) \\
&\label{two}\text{If } E\cong F \text{ it follows } c_i(E)=c_i(F)\\
&\label{three}c_k(E\oplus F)=\sum_{i+j=k}c_i(E)c_j(F).
\end{align}
\end{definition}

Note: if $\phi:\H^*\rightarrow \H^*$ is a functorial endomorphism of
$\H^*$ which is a ring-homomorphism and $c$ is a theory of
characteristic classes, it follows the assignment 
$E\rightarrow \overline{c}_i(E)=\phi(c_i(E))$ is a theory of
characteristic classes.

\begin{example} Let $k\in \Z$ and let $\psi^k_H$ be the ring-endomorphism of $\H^{even}$ defined
by $\psi^k_H(x)=k^rx$ where $x\in \H^{2r}(B)$. Given a theory $c_i(E)$
satisfying Definition \ref{char}
it follows $\overline{c}_i(E)=\psi^k_H(c_i(E))$ is a theory satisfying
Definition \ref{char}.
\end{example}

Note furthermore: Assume $\gamma_1$ is the tautological linebundle on $\P^1$. 
Since we do not assume $c_1(\gamma_1)=z$ where $z$ is the canonical
generator of $\H^2(\P^1,\Z)$ it does not follow that an assignment
$E\rightarrow c_i(E)$ is uniquely determined by the axioms
$\ref{one}-\ref{three}$. We shall see later that the axioms
$\ref{one}-\ref{three}$ is enough to define a Chern-character. 

\begin{theorem}Assume the theory $\H^*$ satisfy the projective bundle
  property. It follows $\H^*$ has a theory of characteristic classes.
\end{theorem}
\begin{proof}We verify the axioms for a theory of characteristic classes.
Axiom \ref{one}: Assume we have a map of rank $n$ bundles $f:F\rightarrow E$
over a map of topological spaces $g:B'\rightarrow B$. We pull back the
equation
\[ u_E^n-c_1(E)u_E^{n-1}+\cdots +(-1)^nc_n(E)=0 \]
in $\H^{2n}(\P(E))$ to get an equation
\[ u_F^n-f^*c_1(E)u_F^{n-1}+\cdots +(-1)^nf^*c_n(E)=0 \]
and by unicity we get $f^*c_i(E)=c_i(F)$. It follows $c_i(E)=c_i(F)$
for isomorphic bundles $E$ and $F$, hence Axiom \ref{two} is ok.
Axiom \ref{three}: Assume $E\cong \oplus_{i=1}^n L_i$ is a
decomposition into linebundles. There is an equation $\prod
(u_E-u_{L_i})$ hence we get a polynomial relation
\[ u_E^n-s_1(u_{L_i})u_E^{n-1}+\cdots +(-1)^ns_n(u_{L_i})=0 \]
in $\H^{2n}(\P(E))$. Since $c_1(L_i)=-u_{L_i}$ it follows
\[ \prod(c(L_i)) =\prod (1+c_1(L_i))=c(E) \]
and this is ok.
\end{proof}

Given a compact topological space $B$. We may consider the
Grothendieck-ring $K_\c^*(B)$ of complex finite-dimensional
vectorbundles. It is defined as the free abelian group on
isomorphism-classes $[E]$ where $E$ is a complex vectorbundle, modulo
the subgroup generated by elements of the type $[E\oplus F
]-[E]-[F]$. It has direct sum as additive operation and tensor product
as multiplication. Assume $E$ is a complex vectorbundle of rank $n$ and let 
\[ \pi: \P(E)\rightarrow B \]
be the associated projective bundle. We have a \emph{projective bundle
  theorem} for complex K-theory:

\begin{theorem} The group $\K^*(\P(E))$ is a free $\K^*(B)$ module of
  finite rank with generator $u$ - the euler class of the tautological line-bundle.
The elements $\{1,u,u^2,..,u^{n-1}\}$ is a free basis.
\end{theorem}
\begin{proof} See Theorem IV.2.16 in \cite{karoubi}.
\end{proof}

As in the case of singular cohomology, we may define characteristic
classes for complex bundles with values in complex K-theory using the
projective bundle theorem: The element $u^n$ satisfies an equation
\[ u^n-c_1(E)u^{n-1}+c_2(E)u^{n-2}+\cdots
+(-1)^{n-1}c_{n-1}(E)u+(-1)^nc_n(E)=0 \]
in $\K^*(\P(E))$. One verifies the axioms defined above are satisfied,
hence one gets characteristic classes $c_i(E)\in
\K^*_\c(B)$ for all $i=0,..,n$

\begin{theorem} The characteristic classes $c_i(E)$ satisfy the
  following properties:
\begin{align}
&f^*c_i(E)=c_i(f^*E) \\
&c_k(E\oplus F)=\sum_{i+j=k}c_i(E)c_j(F)\\
&c_1(L)=1-L
&c_i(L)=0, i>1
\end{align}
where $E$ is any vectorbundle, and $L$ is a line bundle.
\end{theorem}
\begin{proof}See Theorem IV.2.17 in \cite{karoubi}.
\end{proof}

\section{Adams operations and Newton polynomials}

We introduce some cohomology operations in complex K-theory and
Newton-polynomials and prove elementary properties following the book \cite{karoubi}.

Let $\Phi(B)$ be the abelian monoid of elements of the type $\sum
n_i[E_i]$ with $n_i \geq 0$. Consider the bundle
$\lambda^i(E)=\wedge^i E$ and the association
\[ \lambda _t(E)=\sum_{i \geq 0}\lambda^i(E)t^i \]
giving a map
\[ \lambda_t=\Phi(X)\rightarrow 1+t\K_\c^*(B)[[t]]  \]
One checks 
\[ \lambda_t(E\oplus F)=\lambda_t(E)\lambda_t(F) \]
hence the map $\lambda_t$ is a map of abelian monoids, hence gives
rise to a map
\[ \lambda_t: \K_\c^*(B)\rightarrow 1+t\K_\c^*(B)[[t]] \]
from the additive abelian group $\K_\c^*(B)$ to the set of powerseries
with constant term equal to one. Explicitly the map is as follows:
\[ \lambda_t(n[E]-m[F])=\lambda_t(E)^n\lambda_t(F)^{-m}. \]
When $n$ denotes the trivial bundle of rank $n$ we get the explicit
formula
\[ \lambda_t([E]-n)=\lambda_t(E)(1+t)^{-n}.\]
Let $u=t/1-t$. We may define the new powerseries
\[ \gamma_t(E)=\lambda_u(E)=\sum_{k\geq 0}\lambda^i(E)u^i .\]
It follows
\[ \gamma_t(E\oplus F)=\lambda_u(E\oplus
F)=\lambda_u(E)\lambda_u(F)=\gamma_t(E)\gamma_t(F). \]
We may write formally
\[ \gamma_t(E)=\sum_{k\geq 0}\gamma^i(E)t^i \in K^*_\c(B)[[t]] \]
hence it follows that
\[ \gamma^k(E)=\sum_{i+j=k}\gamma^i(E)\gamma^j(E).\]
We get operations
\[ \gamma^i:\K^*_\c(B)\rightarrow \K^*_\c(B) \]
for all $i\geq 1$.
We next define \emph{Newton polynomials} using the elementary
symmetric functions. Let $u_1,u_2,u_3,..$ be independent variables
over the integers $\Z$, and let $Q_k=u_1^k+u_2^k+\cdots +u_k^k$ for 
$k\geq 1$. It follows $Q_k$ is invariant under permutations of the
variables $u_i$: for any $\sigma \in S_k$ we have $\sigma Q_k=Q_k$
hence we may express $Q_k$ as a polynomial in the elementary symmetric
functions
$\sigma_i$:
\[ Q_k=Q_k(\sigma_1,\sigma_2,..,\sigma_k) .\]
We define
\[ s_k(\sigma )=Q_k(\sigma_1,\sigma_2,..,\sigma_k) \]
to be the $k'th$ Newton polynomial in the variables
$\sigma_1,\sigma_2,..,\sigma_k$
where $\sigma_i$ is the $i'th$ elementary symmetric function. 
One checks the following:
\[ s_1(\sigma_1)=\sigma_1, \]
\[ s_2(\sigma_1,\sigma_2)=\sigma_1^2-2\sigma_2 ,\]
and
\[
s_2(\sigma_1,\sigma_2,\sigma_3)=\sigma_1^3-3\sigma_1\sigma_2+3\sigma_3 \]
and so on.

Let $n\geq 1$ and consider the polynomial
\[ p(1)=(1+tu_1)(1+tu_2)\cdots
(1+tu_n)=t^n\sigma_n+t^{n-1}\sigma_{n-1}+\cdots +t\sigma_1+1 \]
where
\[ \sigma_i=\sigma_i(u_1,..,u_n) \]
is the ith elementary symmetric polynomial in the variables 
$u_1,u_2,..,u_n$.

\begin{lemma}\label{elem}  There is an equality
\[Q_k(\sigma_1(u_1,..,u_n),\sigma_2(u_1,..,u_n),..,\sigma_k(u_1,..,u_n))=u_1^k+u_2^k+\cdots
+u_n^k.\]
\end{lemma}
\begin{proof} Trivial.
\end{proof}

Assume we have virtual elements $x=E-n=\oplus^n (L_i-1)$ and
$y=F-p=\oplus^p(R_j-1)$ in complex K-theory $\K^*_\c(B)$. 
We seek to define a cohomology-operation $c$ on complex K-theory using a formal
powerseries
\[ f(u)=a_1u+a_2u^2+a_3u^3+\cdots \in \Z[[u]]. \]
We define the element
\[
c(x)=a_1Q_1(\gamma^1(x))+a_2Q_2(\gamma^1(x),\gamma^2(x))+a_3Q_3(\gamma^1(x),\gamma^2(x),\gamma^3(x))+\cdots
\]

\begin{proposition} Let $L$ be a linebundle. Then
  $\gamma_t(L-1)=1+t(L-1)=1-c_1(L)t$. Hence $\gamma^1(L-1)=L-1$ and
  $\gamma^i(L-1)=0$ for $i>1$. 
\end{proposition}
\begin{proof} We have by definition
\[ \gamma_t(E)=\lambda_u(E)=\sum_{k\geq 0}\lambda^k (E)u^k=\sum_{k\geq
  0}\lambda^k(E)(t/1-t)^k. \]
We have that 
\[ \gamma_t(nE-mF)=\lambda_u(E)^{n}\lambda_u(F)^{-m}. \]
We get 
\[ \gamma_t(L-1)=\lambda_u(L)\lambda_u(1)^{-1} .\]
We have
\[ \lambda_t(n)=(1+t)^n \]
hence
\[ \gamma_t(n)=\lambda_u(n)=(1+u)^n=(1+t/1-t)^n=(1-t)^{-n}.\]
We get:
\[\gamma_t(L-1)=\gamma_t(L)\gamma_t(1)^{-1}=\lambda_u(L)(1-t)^{-1}=\]
\[(1+Lu)(1-t)^{-1}=(1+L(t/t-1))(1-t)^{-1}= \]
\[ \frac{1+t(L-1)}{1-t}(1-t)=1+t(L-1)=1-c_1(L)t.\]
And the proposition follows.
\end{proof}
Note: if $x=L-1$ we get
\[ c(x)=\sum_{k\geq 0}a_kQ_k(\gamma^1(x),\gamma^2(x),..,\gamma^k(x))=
\]
\[ \sum_{k\geq 1}a_k Q_k(\gamma^1(x),0,...,0)=\sum_{k\geq
  1}a_k\gamma^1(x)^k= \]
\[ \sum_{k\geq 1}a_k(L-1)^k=\sum_{k\geq 0}(-1)^ka_kc_1(L)^k.\]

We state a Theorem:

\begin{theorem}\label{proj_ktheory} Let $E\rightarrow B$ be a complex vectorbundle on a
  compact topological space $B$. There is a map $\pi:B'\rightarrow B$
  such that $\pi^*E$ decompose into linebundles, and the map
  $\pi^*:\H^*(B)\rightarrow \H^*(B')$ is injective.
\end{theorem}
\begin{proof} See \cite{karoubi} Theorem IV.2.15.
\end{proof}

Note: By \cite{karoubi} Proposition II.1.29 there is a split exact sequence
\[ 0\rightarrow K'_\c(B)\rightarrow \K^*_\c(B)\rightarrow
\H^0(B,\Z)\rightarrow 0 \]
hence the group $K'_\c(B)$ is generated by elements of the form $E-n$
where $E$ is a rank $n$ complex vectorbundle.

\begin{proposition} The operation $c$ is additive, i.e for any $x,y\in
  \K^*_\c(B)$ we have
\[ c(x+y)=c(x)+c(y). \]
\end{proposition}
\begin{proof} The proof follows the proof in \cite{karoubi},
  Proposition IV.7.11. We may by the remark above assume $x=E-n$ and
  $y=F-p$ where $x,y\in K'_\c(B)$. We may also from Theorem \ref{proj_ktheory} assume $E=\oplus^n L_i$
  and $F=\oplus^p R_j$ where $L_i,R_j$ are linebundles. We get the
  following:
\[ \gamma_t(x+y)=\prod^n \gamma_t(L_i-1)\prod^p
\gamma_t(R_j-1)=\prod^n(1+tu_i)\prod^p (1+tv_j)=\]
\[t^{n+p}\sigma_{n+p}(u_1,..,u_n,v_1,..,v_p)+t^{n+p-1}\sigma_{n+p-1}(u_1,..,u_n,v_1,..,v_p)+
\] 
\[ \cdots +t\sigma_1(u_1,..,u_n,v_1,..,v_p)+1 \]
hence
\[ \gamma^i(x+y)=\sigma_i(u_1,..,u_n,v_1,..,v_p).\]
We get:
\[ Q_k(\gamma^1(x+y),..,\gamma^k(x+y))=Q_k(\sigma_1(u_i,v_j),..,\sigma_k(u_i,v_j)) \]
which by Lemma \ref{elem} equals
\[u_1^k+\cdots u_n^k+v_1^k+\cdots
v_p^k=Q_k(\sigma_1(u_i),..,\sigma_k(u_i))+Q_k(\sigma_1(v_j),..,\sigma_k(v_j))=\]
\[ Q_k(\gamma^i(x))+Q_k(\gamma^i(y)). \]
We get:
\[c(x+y)=\sum_{k\geq 0}a_k Q_k(\gamma^i(x+y))= \]
\[ \sum_{k\geq 0} a_kQ_k(\gamma^i(x))+\sum_{k\geq 0}a_kQ_k(\gamma^i(y))=c(x)+c(y) \]
and the claim follows.

\end{proof}

We may give an explicit and elementary construction of the Adams-operations:

\begin{theorem} Let $k\geq 1$. There are functorial operations
\[ \psi^k:\K^*_\c(B)\rightarrow \K^*_\c(B) \]
with the properties
\begin{align}
& \psi^k(x+y)=\psi^k(x)+\psi^k(y) \\
& \psi^k(L)=L^k \\
& \psi^k(xy)=\psi^k(x)\psi^k(y)\\ 
& \psi^k(1)=1
\end{align}
where $L$ is a line bundle. The operations $\psi^k$ are the only
operations that are ring-homomorphisms - the \emph{Adams operations}
\end{theorem}
\begin{proof} We need: 
\[ \psi^k(L-1)=\psi^k(L)-\psi^k(1)=L^k-1. \]
We have in K-theory:
\[ L^k-1=(L-1+1)^k-1=\sum_{i\geq 0}\binom{k}{i}(L-1)^{k-i}{1}^i-1=      \]
\[ \binom{k}{1}(L-1)+\binom{k}{2}(L-1)^2+\cdots +\binom{k}{k}(L-1)^k. \]
We get the series
\[ c=\sum_{i=1}^k\binom{k}{i}u^k \in \Z[[u]] .\]
The following operator
\[ \psi^k=\sum_{i=1}^k \binom{k}{i}Q_i(\gamma^1,...,\gamma^i) \]
is an explicit construction of the Adams-operator. One may verify the
properties in the theorem, and the claim follows.
\end{proof}

Assume $E,F$ are complex vectorbundles on $B$ and consider the
Chern-polynomial
\[ c_t(E\oplus F)=1+c_1(E\oplus F)t+\cdots +c_N(E\oplus F)t^N \]
where $N=rk(E)+rk(F)$. Assume there is a decomposition $E=\oplus^n
L_i$ and $F=\oplus^p R_j$ into linebundles. We get a decomposition
\[ c_t(E\oplus F)=\prod c_t(L_i)\prod c_t(R_j)=(1+a_1t)(1+a_2t)\cdots
(1+b_1t)\cdots (1+b_pt) \]
where $a_i=c_1(L_i), b_j=c_1(R_j)$. We get thus
\[ c_i(E\oplus F)=\sigma_i(a_1,..,a_n,b_1,..,b_p).\]
Let 
\[Q_k=u_1^k+\cdots +u_k^k=Q_k(\sigma_1,..,\sigma_k) \]
where $\sigma_i$ is the ith elementary symmetric function in the $u_i$'s.

\begin{proposition} \label{qk} The following holds:
\[ Q_k(c_1(E\oplus F),..,c_k(E\oplus F))=Q_K(c_i(E))+Q_k(c_i(F)).\]
\end{proposition}
\begin{proof} We have 
\[Q_k(c_i(E\oplus F))=Q_k(\sigma_i(a_i,b_j))= \]
\[a_1^k+\cdots a_n^k+ b_1^k+\cdots b_p^k =Q_k(c_i(E))+Q_k(c_i(F)) \]
and the claim follows.
\end{proof}

\section{The Chern-character and cohomology operations}

We construct a Chern-character with values in singular cohomology,
using Newton-polynomials and characteristic classes following \cite{karoubi}. 
The $k'th$ Newton-classe $s_k(E)$ of a complex vectorbundle will be
defined using characteristic classes of $E$: $c_1(E),..,c_k(E)$ and the
$k'th$ Newton-polynomial $s_k(\sigma_1,..,\sigma_k)$. We us this
construction to define the Chern-character $Ch(E)$ of the vectorbundle
$E$. 

We first define \emph{Newton polynomials} using the elementary
symmetric functions. Let $u_1,u_2,u_3,..$ be independent variables
over the integers $\Z$, and let $Q_k=u_1^k+u_2^k+\cdots +u_k^k$ for 
$k\geq 1$. It follows $Q_k$ is invariant under permutations of the
variables $u_i$: for any $\sigma \in S_k$ we have $\sigma Q_k=Q_k$
hence we may express $Q_k$ as a polynomial in the elementary symmetric
functions
$\sigma_i$:
\[ Q_k=Q_k(\sigma_1,\sigma_2,..,\sigma_k) .\]
We define
\[ s_k(\sigma )=Q_k(\sigma_1,\sigma_2,..,\sigma_k) \]
to be the $k'th$ Newton polynomial in the variables
$\sigma_1,\sigma_2,..,\sigma_k$
where $\sigma_i$ is the $i'th$ elementary symmetric function. 
One checks the following:
\[ s_1(\sigma_1)=\sigma_1, \]
\[ s_2(\sigma_1,\sigma_2)=\sigma_1^2-2\sigma_2 ,\]
and
\[
s_2(\sigma_1,\sigma_2,\sigma_3)=\sigma_1^3-3\sigma_1\sigma_2+3\sigma_3 \]
and so on.

Assume we have a cohomology theory $\H^*$ satisfying the projective
bundle property. One gets characteristic 
classes $c_i(E)$ for a complex vectorbundle $E$ on $B$: 
\[ c_i(E)\in \H^{2i}(B) .\]
Let the class $S_k(E)=s_k(c_1(E),c_2(E),..,c_k(E))\in \H^{2k}(B)$ be the
$k'th$ Newton-class of the bundle $E$. 
One gets:
\[ s_k(\sigma_1,0,..,0)=\sigma_1^k \]
for all $k\geq 1$. Assume $E,F$ linebundles. We see that 
\[ S_2(E\oplus F)=c_1(E \oplus F)^2-2c_2(E\oplus F)= \]
\[ (c_1(E)+c_1(F))^2-2(c_2(E)+c_1(E)c_1(F)+c_2(F))=\]
\[ c_1(E)^2+2c_1(E)c_1(F)+c_1(F)^2-2c_2(E)-2c_1(E)c_1(F)-2c_2(F)=\]
\[ c_1(E)^2-2c_2(E)+c_1(F)^2-2c_2(F)=S_2(E)+S_2(F).\]
This holds in general:
\begin{proposition} \label{newton1} For any vectorbundles $E,F$ we have the formula
\[ S_k(E\oplus F)=S_k(E)+S_k(F).\]
\end{proposition}
\begin{proof} This follows from \ref{qk}.
\end{proof}

Let $\K_\c^*(B)$ be the Grothendieck-group of complex vectorbundles on $B$,
i.e the free abelian group modulo exact sequences $K_\c^*(B)=\oplus \Z
[E]/U$ where $U$ is the subgroup generated by elements $[E\oplus
F]-[E]-[F]$.

\begin{definition} The class
\[ Ch(E)=\sum_{k\geq 0}\frac{1}{k!}S_k(E)\in \H^{even}(B) \]
is the \emph{Chern-character} of $E$.
\end{definition}

\begin{lemma} The Chern-character defines a group-homomorphism
\[ Ch: \K_\c^*(B)\rightarrow \H^{even}(B) \]
between the Grothendieck group $\K^*_\c(B)$ and the even cohomology of
$B$ with rational coefficients.
\end{lemma}
\begin{proof} By Proposition \ref{newton1} we get the following: For
  any $E,F$ we have
\[Ch(E\oplus F)=\sum_{k\geq 0}\frac{1}{k!}s_k(E\oplus F)=\sum_{k\geq
  0}\frac{1}{k!}(s_k(E)+s_k(F))=\]
\[ \sum_{k\geq 0}\frac{1}{k!}s_k(E)+\sum_{k\geq
  0}\frac{1}{k!}s_k(F)=Ch(E)+Ch(F).\]
We get 
\[ Ch([E\oplus F ]-[E]-[F])=Ch(E\oplus F)-Ch(E)-Ch(F)=0 \]
and the Lemma follows.
\end{proof}

\begin{example}
Given a real continuous vectorbundle $F$ on $B$ there exist
Stiefel-Whitney classes $w_i(F)\in \H^i(B,\Z/2)$ (see \cite{mil}) satisfying the
necessary conditions, and we may define a ``Chern-character''
\[ Ch:\K_\r^*(B)\rightarrow \H^*(B,\Z/2) \]
by
\[ Ch(F)=\sum_{k\geq 0}Q_k(w_1(F),..,w_k(F)) .\]
This gives a well-defined homomorphism of abelian groups because of the universal
properties of the Newton-polynomials and the fact $\H^*(B,\Z/2)$ is
commutative. 
The formal properties of the Stiefel-Whitney classes $w_i$ ensures
that for real bundles $E,F$ Proposition \ref{qk} still holds:
We have the formula
\[ Q_k(w_i(E\oplus F))=Q_k(w_i(E))+Q_k(w_i(F)). \]
\end{example}

Since $S_k(\sigma_1,0,...,0)=\sigma_1^k$ we get the following: When
$E,F$ are linebundles we have:
\[ S_k(E\otimes F)=S_k(c_1(E\otimes F),0,..,0)=(c_1(E\otimes
F))^k=(c_1(E)+c_1(F))^k=\]
\[
\sum_{i+j=k}\binom{i+j}{i}c_1(E)^ic_1(F)^j=\sum_{i+j=k}\binom{i+j}{i}S_i(E)S_j(F).\]
This property holds for general $E,F$:

\begin{proposition} \label{newton2} Let $E,F$ be complex vectorbundles on a
  compact topological space $B$. Then the following formulas hold:
\begin{align}
& S_k(E\otimes F)=\sum_{i+j=k}\binom{i+j}{i}S_i(E)S_j(E) 
\end{align}
\end{proposition}
\begin{proof} We prove this using the splitting-principle and Proposition
  \ref{newton1}. Assume $E,F$ are complex vectorbundles on $B$ and
  $f:B'\rightarrow B$ is a map of topological spaces such that
  $f^*E=\oplus_i L_i, f^*F=\oplus_j M_j $ where $L_i,M_j$ are
  linebundles and the pull-back map $f^*:\H^*(B)\rightarrow \H^*(B')$
  is injective. We get the following calculation: 
\[ f^*S_k(E\otimes F)=S_k(f^*E\otimes F)=S_k(\oplus L_i\otimes M_j) \]
hence by Lemma \ref{newton1} we get
\[ \sum_{i,j}S_k(L_i\otimes M_j)=\sum_i(\sum_j S_k(L_i\otimes M_j))
=\]
\[ \sum_i \sum_j
\sum_{u+v=k}\binom{u+v}{u}S_u(L_i)S_v(M_j)=\]
\[ \sum_i \sum_{u+v=k}\binom{u+v}{u}S_u(L_i)S_v(\oplus M_j) =\]
\[ \sum_{u+v=k}\binom{u+v}{u} S_u(\oplus L_i)S_v(\oplus M_j)= \]
\[
\sum_{u+v=k}\binom{u+v}{u}S_u(f^*E)S_v(f^*F)=f^*\sum_{u+v=k} \binom{u+v}{u}S_u(E)S_v(F),
\]
and the result follows since $f^*$ is injective.
\end{proof}

\begin{theorem} The Chern-character defines a ring-homomorphism
\[ Ch: K_\c^*(B)  \rightarrow \H^{even}(B) .\]
\end{theorem}
\begin{proof} From Proposition \ref{newton2} we get:
\[ Ch(E\otimes F)=\sum_{k \geq 0}\frac{1}{k!}S_k(E\otimes F)=\]
\[ \sum_{k\geq 0}\frac{1}{k!}\sum_{i+j=k}\binom{i+j}{i}S_i(E)S_j(F)=\]
\[ (\sum_{k\geq 0}\frac{1}{k!}S_k(E))(\sum_{k\geq
  0}\frac{1}{k!}S_k(F)=Ch(E)Ch(F) \]
and the Theorem is proved.
\end{proof}

\begin{example}
For complex K-theory $\K^*_\c(B)$ we have for any
complex vectorbundle $E$ characteristic classes $c_i(E)\in \K^*_\c(B)$
satisfying the neccessary conditions, hence we get a group-homomorphism
\[ Ch_\Z:\K_\c^*(B)\rightarrow \K^*_\c(B) \]
defined by
\[ Ch_\Z(E)=\sum_{k\geq 0}Q_k(c_1(E),..,c_k(E)) .\]
If we tensor with the rationals, we get a ring-homomorphism
\[ Ch_\Q:\K^*_\c(B)\rightarrow \K^*_\c(B)\otimes \Q \]
defined by
\[ Ch(E)=\sum_{k\geq 0}\frac{1}{k!}Q_k(c_1(E),..,c_k(E)) .\]
\end{example}

\begin{theorem} Let $B$ be a compact topological space. The
Chern-character 
\[ Ch^\Q: K_\c^*(B)\otimes \Q \rightarrow \H^{even}(B,\Q) \]
is an isomorphism. Here $\H^*(B,\Q)$ denotes singular cohomology with
rational coefficients.
\end{theorem}
\begin{proof} See \cite{karoubi}.
\end{proof}

The Chern-character is related to the Adams-operations in the
following sense: There is a ring-homomorphism
\[ \psi^k_H:\H^{even}(B)\rightarrow \H^{even}(B) \]
defined by
\[ \psi^k_H(x)=k^rx \]
when $x\in \H^{2r}(B)$. 
The Chern-character respects these cohomology operations in the following sense:

\begin{theorem}\label{main} There is for all $k\geq 1$ a commutative diagram
\[ 
\diagram  K^*_\c(B) \rto^{Ch}\dto^{\psi^k} & \H^{even}(B)
\dto^{\psi^k_H} \\
K^*_\c(B) \rto^{Ch}  & \H^{even}(B) 
\enddiagram 
\]
where $\psi^k$ is the Adams operation defined in the previous section.
\end{theorem} 
\begin{proof} The proof follows Theorem V.3.27 in \cite{karoubi}: 
We may assume $L$ is a linebundle and we get the
  following calculation: $\psi^k(L)=L^k$ and $c_1(L^k)=kc_1(L)$ hence
\[ Ch(\psi^k(L))=exp(kc_1(L))=\sum_{i\geq 0}\frac{1}{i!}k^ic_1(L)^i=
\]
\[ \psi^k_H(exp(c_1(L)))=\psi^k_H(Ch(L) \]
and the claim follows.
\end{proof}

Hence the Chern-character is a morphism of cohomology-theories respecting
the additional structure given by the Adams and
Steenrod-operations.

\end{document}